\documentclass[12pt,a4paper]{article}
\usepackage[T1]{fontenc}
\usepackage{mathptmx}
\usepackage{courier}
\usepackage[dvips]{graphicx}
\usepackage{psfrag}
\usepackage{amsmath,amssymb}        
\usepackage{latexsym}
\usepackage{url}

\def\ds{\displaystyle}
\def\e{\epsilon}

\def\E{\ensuremath{E}}

\newtheorem{theorem}{Theorem}
\def\ds{\displaystyle}
\def\e{\epsilon}
\def\FBF{free boundary formulation}

\def\www{.9\textwidth}

\title{Free Boundary Formulation for Boundary Value Problems on Semi-Infinite Intervals:\\
 An up to Date Review}

\author{Riccardo Fazio \\
Department of Mathematics, Computer Science\\ Physical Sciences and Earth Sciences,\\
University of Messina \\
Viale F. Stagno D'Alcontres, 31 \\
98166 Messina, Italy \\
E-mail: rfazio@unime.it \\
Home-page: http://mat521.unime.it/fazio}
\date{\today}
\begin{document}
\maketitle
\begin{abstract}
In this paper, we propose a review of the \FBF \ for BVPs defined on semi-infinite intervals.
The main idea and theorem are illustrated, for the reader convenience, by using a class of second-order BVPs.
Moreover, we are able to show the effectiveness of the proposed approach using two examples where the exact solution both for the BVPs and their \FBF \ are available. 
Then, we describe the \FBF \ for a general class of BVPs governed by an $n$-order differential equation. 
In this context, we report three problems solved using the \FBF .
The reported numerical results, obtained by the iterative transformation method or the Keller's second-order finite difference method, are found to be in very good agreement with those available in the literature.
The last result of this research is that, in order to orient the interested reader, we provide an extensive bibliography.
Of course, we may aspect further and more interesting applications of the \FBF \ in the future.
\end{abstract}

\bigskip


\section{Introduction and bibliography}
Usually, when dealing for the first time with a boundary value problem (BVP) defined on a semi-infinite interval, the applied scientist does not know the exact or even an approximate solution.
As a consequence, he often is tempted to try for a numerical solution to the problem.
Therefore, along the years' several approaches have been proposed in order to solve such a problem numerically.

The oldest and simplest approach, to deal with the problems we are facing, is to replace the original problem by one defined on a finite interval, where a finite value, the so-called truncated boundary, is used instead of infinity. 
This approach was used, for instance, to get the tabulated numerical solution \cite{Horwarth:1938:SLB} and \cite[p. 136]{Goldstein:1938:MDF} of the Blasius problem \cite{Blasius:1908:GFK}.
However, to get an accurate solution a comparison of numerical results obtained for several values of the truncated boundary is necessary as suggested by Fox \cite[p. 92]{Fox} and by Collatz \cite[pp. 150-151]{Collatz}.
Moreover, in some cases, accurate solutions can be found only by using very large values of the truncated boundary.
This is, for instance, the case for the branches of the von Karman swirling flows where values of truncated boundaries up to several hundred were used by Lentini and Keller \cite{Lentini:KSF:1980}.

The limitation of the above classical approach has lead some researchers de Hoog and Weiss \cite{deHoog:1980:ATB}, Lentini and Keller \cite{Lentini:BVP:1980} and
Markowich \cite{Markowich:TAS:1982,Markowich:ABV:1983} to develop a theory for defining the asymptotic boundary conditions to be imposed at a given value of the truncated boundary. 
Those asymptotic boundary conditions are derived by a preliminary asymptotic analysis involving the Jacobian matrix of the right-hand side of the governing equations evaluated at infinity.
The main idea of this asymptotic boundary conditions approach is to project the solution into the manifold of bounded solutions.
By using the same value of the truncated boundary, a more accurate numerical solution can be found by this approach than that obtained by the classical approach, because in the first case the imposed boundary conditions are obtained from the asymptotic behaviour of the solution.
However, we should note that this approach is not straightforward, see the remarks by Ockendon \cite{Ockendon}, and that for nonlinear problems highly nonlinear asymptotic boundary conditions usually result.
Asymptotic boundary conditions have been applied successfully to the numerical approximation of the so-called \lq \lq connecting orbits\rq \rq \ problems of dynamical systems, see Beyn \cite{Beyn:1990:GBN,Beyn:1990:NCC,Beyn:1992:NMD}.
Those problems are of interest, not only in connection with dynamical systems, but also in the study of traveling wave solutions of partial differential equations of parabolic and hyperbolic type as shown by Beyn \cite{Beyn:1990:NCC}, Friedman \cite{Friedman:1991:NCC}, Bai et al. \cite{Bai:1993:NCH}, and Liu et al. \cite{Liu:1997:CCH}.

A different approach, for the numerical solution of BVPs defined on a semi-infinite, is to consider a free boundary formulation of the given problem, where the unknown free boundary can be identified with a truncated boundary. 
In this approach, the free boundary is unknown and has to be found as part of the solution. 
This free boundary approach overcomes the need for an  \textit{a priori} definition of the truncated boundary. 
Free BVPs represent a numerical challenge because they are always non-linear as pointed out first by Landau \cite{Landau:1950:HCM}.
However, a \FBF \ has been successfully applied to several problems in the applied sciences: namely,  the Blasius problem by Fazio \cite{Fazio:1992:BPF}, a two-dimensional stagnation point flow by Ariel \cite{Ariel:FEF:1993}, the Falkner-Skan model by Fazio \cite{Fazio:1994:FSEb}, by Zhang and Cheng \cite{Zhang:IMS:2009} and by Zhu et al. \cite{Zhu:NSF:2009}, and the model describing a fluid flowing around a slender parabola of revolution by Fazio \cite{Fazio:1996:NAN} in boundary layer theory, the computation of a two-dimensional homoclinic connecting orbit by Fazio \cite{Fazio:2002:SFB}, and a problem related to the deflection of a semi-infinite pile embedded in soft soil by Fazio \cite{Fazio:2003:FBA}.
The last problem is of interest in foundation engineering, for instance, in the design of drilling rigs above the ocean floor, see Lentini and Keller \cite{Lentini:BVP:1980} and the references quoted therein. 

A different way to avoid the definition of a truncated boundary is to apply coordinate transforms.
Coordinate transforms have been applied successfully to the numerical solution of ordinary and partial differential equations on unbounded domains, see Grosch and Orszag \cite{Grosch:NSP:1977}, Koleva\cite{Koleva:NSH:2006} or Fazio and Jannelli \cite{Fazio:2014:FDS}.

The main idea and theorem related to the \FBF \ are illustrated using a class of second-order BVPs.
In this context, we show in full details the application of the \FBF \ to two example of BVPs defined on semi-infinite intervals.
In both cases, we are able to provide the exact solution of both the BVP and its \FBF .
Therefore, these problems can be used as benchmarks for the numerical methods applied to BVPs on a semi-infinite intervals and to free BVPs. 
In this context, sometimes, it is possible to solve a given free BVP non-iteratively, see the survey by Fazio \cite{Fazio:1998:SAN}, whereas BVPs are usually solved iteratively.
Here, for two classes of free BVPs, we define non-iterative initial value methods which, in the literature, are referred to as non-iterative transformation methods (ITMs).
Indeed, non-ITMs can be defined within Lie's group invariance theory.
For the group invariance theory, the interested reader is referred to Bluman and Cole \cite{Bluman:1974:SMD}, Bluman and Kumei \cite{Bluman:1989:SDE}, Barenblatt \cite{Barenblatt:1996:SSI}, or Dresner \cite{Dresner:1999:ALT}.

A review paper \cite{Fazio:1998:SAN} by this authors proposed a brief description related to the main topic of this survey, the numerical example there was concerning with a fluid flow around a slender parabola of revolution mimic an airplane engine.

The main goal of this paper is to provide evidence of the effectiveness of the \FBF for the class of BVPs defined on semi-infinite intervals defined by
\begin{align}\label{eq:class1}
& {\ds \frac{d^nu}{dx^n}} + f\left(x, u, \cdots, {\ds \frac{d^{n-1}u}{dx^{n-1}}}\right) = 0 \nonumber \\
& {\ds \frac{d^ku}{dx^k}} (0) = u_k \ , \qquad \mbox{for} \qquad k = 0, 1, \cdots, n-2 \\
& {\ds \frac{d^ru}{dx^r}} (\infty) = u_{\infty} \ , \nonumber
\end{align}
where $n$ is a positive integer bigger than one, $f(\cdot,\cdots,\cdot)$ is a given function of its arguments, $u_k$, for $k = 0, 1, \cdots, n-2$, $r \in \{0, 1, \cdots, n-1\}$, and $u_\infty$ are given constants.
However, to let the interested reader gain confidence with the \FBF we are going to start, in the next section, with the simplest problem of this kind and its \FBF .

\section{Free boundary formulation main idea}
In order to explain the main idea behind our \FBF , we consider the simplest sub-class of BVPs that belongs to (\ref{eq:class1}), namely 
\begin{align}\label{eq:BVPs}
&{\ds \frac{d^2u}{dx^2}
+ g\left(x, u,\frac{du}{dx}\right)} = 0 \ , \qquad x \in [0, \infty) \nonumber \\[-1.5ex]
&\\[-1.5ex]
&{\ds u(0) = u_0 \ ,  \qquad u(\infty)}
= u_\infty \nonumber
\end{align}
where $g(\cdot,\cdot,\cdot)$ is a given function of its arguments, and $u_0$ and $u_\infty$ are given constants. 
If we can assume that the first derivative of $u(x)$ goes monotonically to zero at infinity, then we replace the problem (\ref{eq:BVPs}) with its \FBF 
\begin{align}\label{eq:BVPs:FBF}
&{\ds \frac{d^2u_\e}{dx^2}
+ g\left(x, u_\e,\frac{du_\e}{dx}\right)} = 0 \ , \qquad x \in [0, x_\e] \nonumber \\[-1.5ex]
&\\[-1.5ex]
&{\ds u_\e(0) = u_0 \ ,  \qquad u_\e(x_\e)}
= u_\infty \ , \quad \frac{du_\e}{dx} (x_\e) = \e \nonumber
\end{align}
where $x_\e$ is an unknown free boundary and $0 \le |\e| \ll 1$ is a parameter.

We have to remark here that monotonic properties of the solution, its first and second derivative have been demonstrated by Countyman and Kannan \cite{Countryman:1994:CNB}, for the class of problems in (\ref{eq:BVPs}) where $g$ depends exclusively on $u$.    

The following theorem provides, under suitable smoothness conditions, the order of convergence (and the uniform convergence) of the solution of
(\ref{eq:BVPs:FBF}) to the solution of (\ref{eq:BVPs}).

\medskip
\noindent
\begin{theorem}\label{Th:Conv} 
Suppose $ u_{\epsilon}(x) $ and $ {\frac{\partial u_{\epsilon}}{\partial \epsilon }} (x) $ are continuous functions with respect to $ \epsilon $ (and also with respect to $ x $ in the related free boundary domain $ [0, x_\epsilon ] $) and that $ | \epsilon_1 | < | \epsilon_2 | \Rightarrow [0, x_{\epsilon _2}] \subset
[0, x_{\epsilon _1}] $ at least in a non-empty interval including $ \epsilon = 0 $, then 
\begin{eqnarray*}
& || u_{\e} (x) - u(x) || \leq L | \e |
\end{eqnarray*}
where $ || \cdot || $ is the maximum norm on $[0, x_\e]$ and $L$ is a positive constant independent on $ \e $.
\end{theorem}

\medskip
\noindent
The proof of this Theorem can be obtained along the lines of the proof for the convergence Theorem stated in Fazio \cite{Fazio:1996:NAN} for a \FBF \ for a class of problems governed by a third-order differential equation.

The \FBF \ allows us to embed a BVP in (\ref{eq:BVPs}) into a class of problems involving the control parameter $\e$.
When we solve the \FBF \ (\ref{eq:BVPs:FBF}) numerically, we can fix a very small value of $|\e|$ and apply a grid refinement to verify whether the numerical results agree within a prefixed number of significant digits. 
Also, it is possible to fix a step size and let $\e$ goes to zero and verify whether $u_\e(x) \rightarrow u(x)$ together with $x_\e \rightarrow \infty$. 
Usually, it suffices to take $|\e| \in \left\{10^{-1}\right.$, $10^{-2}$, $10^{-3}$, $10^{-4}$,
 $10^{-5}$, $10^{-6}$, $\left. \dots \right\}$
and compare the obtained numerical results.
Let us remark here that sometimes it is possible to solve the \FBF \ non-iteratively, see the survey by Fazio \cite{Fazio:1998:SAN}, whereas the BVP (\ref {eq:BVPs}) is usually solved iteratively.
 
\section{Two examples for the \FBF}
\begin{figure*}[!hbt]
\centering
\psfrag{x}[][]{$x$} 
\psfrag{u}[][]{$u(x)$} 
\psfrag{p01}[][]{$P = 0.1$} 
\psfrag{p1}[][]{$P = 1$} 
\psfrag{p10}[][]{$P = 10$} 
\includegraphics[width=\www]{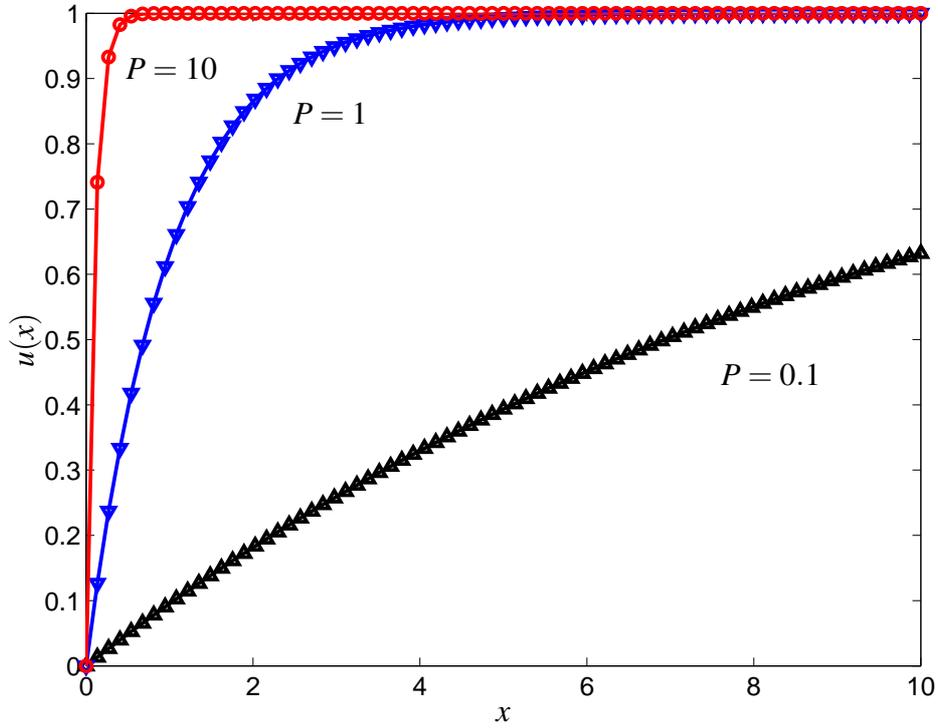}
\caption{The solution (\ref{eq:sol1}) for different values of  $P$. The symbols stand for:
$\circ$ $P=10$, $\triangledown$ $P=1$, and $\vartriangle$ $P=0.1$.}
\label{fig:BVPes1}
\end{figure*}

Let us consider, now, a first example of a BVP defined on an semi-infinite interval.
So, we consider the linear problem
\begin{align}\label{eq:ex1}
& {\ds \frac{d^2u}{dx^2} + P \frac{du}{dx}} = 0 \ , \qquad x \in [0, \infty) \nonumber \\[-1.5ex]
& \\[-1.5ex]
& u(0) = 0 \ , \qquad u(\infty) = 1 \nonumber 
\end{align}
where $P$ is a positive constant.
The solution of (\ref{eq:ex1}) is easily found to be
\begin{equation}\label{eq:sol1}
u(x) = 1-e^{-Px}
\end{equation}
so that the missing initial condition is equal to $P$, that is $\frac{du}{dx}(0) = P$.
Figure \ref{fig:BVPes1} shows the solution (\ref{eq:sol1}) of the BVP (\ref{eq:ex1}) for different values of  $P$.
The bigger is the value of $P$, the harder is to solve the BVP numerically.
In fact, for large values of $P$, the solution has a fast transient for small values of $x$. 

\begin{figure*}[!hbt]
\centering
\psfrag{x}[][]{$x$} 
\psfrag{u}[][]{$u(x), u_\e(x)$} 
\psfrag{P}[][]{$P = 1$} 
\includegraphics[width=\www]{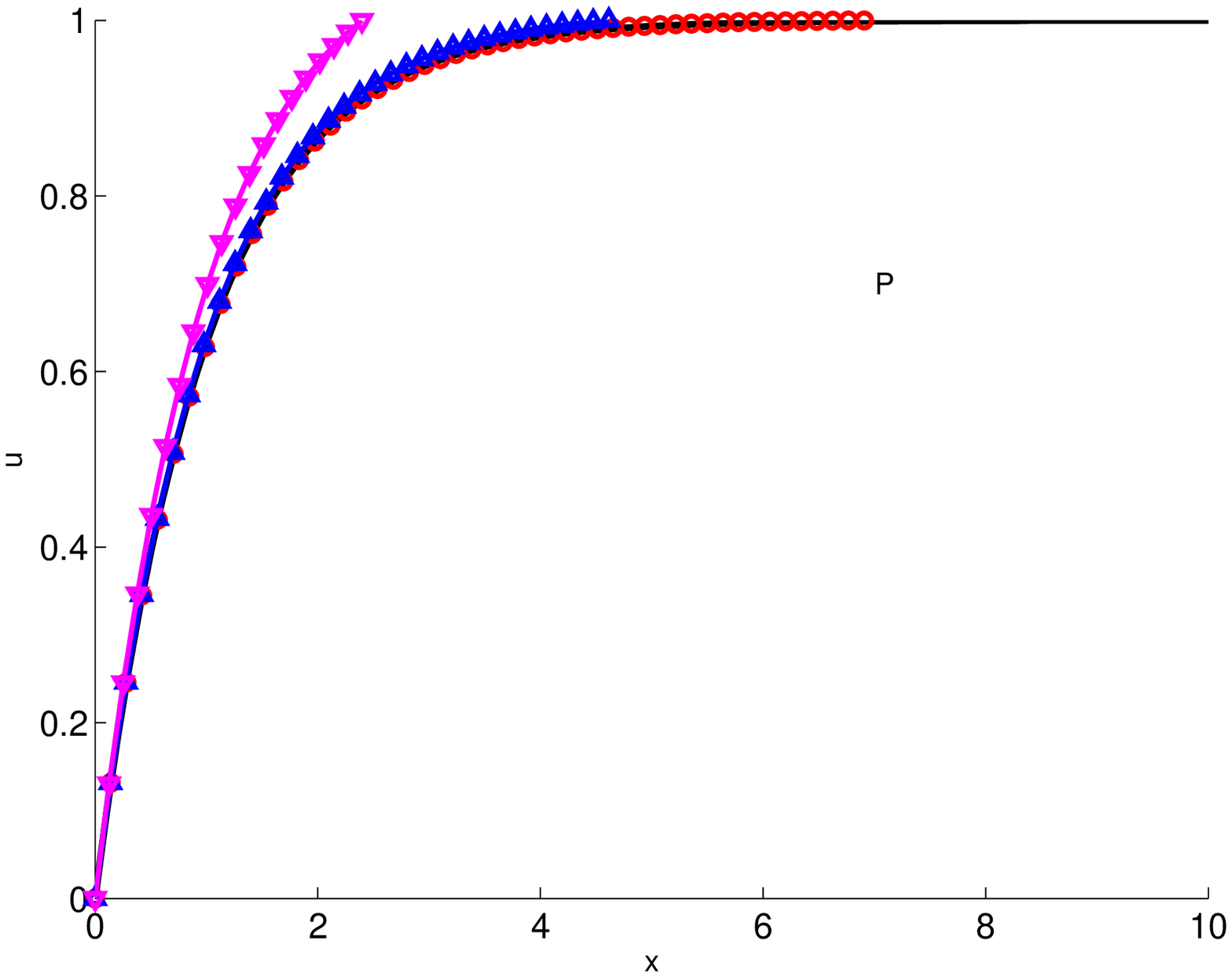}
\caption{The solution (\ref{eq:solfbf1}) for $P = 1$ and different values of  $\e$. The symbols stand for:
$\mathbf{-}$ the exact solution, $\triangledown$, $\vartriangle$  and $\circ$ the free boundary solution $u_\e$ with $\e=0.1$, $\e=0.01$ and $\e=0.001$, respectively.}
\label{fig:FBFes1}
\end{figure*}

Let us consider now the \FBF \ for (\ref{eq:ex1})
\begin{align}\label{eq:ex1:fbf}
& {\ds \frac{d^2u_\e}{dx^2} + P \frac{du_\e}{dx}} = 0 \ , \qquad x \in [0, x_\e] \nonumber \\[-1.5ex]
& \\[-1.5ex]
& u_\e(0) = 0 \ , \qquad u_\e(x_\e) = 1 \ , \qquad \frac{du_\e}{dx}(x_\e) = \e \ , \nonumber 
\end{align}
with $0 \le \e \ll 1$.
The solution of (\ref{eq:ex1:fbf}) is given by
\begin{equation}\label{eq:solfbf1}
u_\e (x) = {\ds \frac{P+\e}{P}\left(1-e^{-Px}\right)} \ , \quad x_\e = - {\ds \frac{1}{P} \ln\left(\frac{\e}{P+\e}\right)} \ .
\end{equation}
Therefore, we can easily verify that as $\e$ goes to zero the solution $u_\e(x)$ of the \FBF \
(\ref{eq:ex1:fbf}) converges to the solution $u(x)$ of the original problem (\ref{eq:ex1}) and the free boundary $x_\e$ goes to infinity.
Moreover, we realize that the obtained approximation becomes the more accurate the more $\e$ is near zero, see figure \ref{fig:FBFes1}.

\noindent
{\bf Remark.} the same exact solutions (\ref{eq:sol1}) and (\ref{eq:solfbf1}) are still valid if we replace the governing differential equation, in the BVP (\ref{eq:ex1}) and its \FBF \ (\ref{eq:ex1:fbf})
with the non-autonomous one
\begin{equation}
{\ds \frac{d^2u}{dx^2} + P^2 e^{-P x}} = 0 \ ,
\end{equation}
where we substitute $u = u_\e$ in the free boundary case.

Replacing a linear problem with a non-linear one can be justified, from a numerical viewpoint, only by considering that in this way we overcome the singularity related to the boundary condition prescribed at infinity.
Of course, when the original problem is a non-linear one a \FBF \ for it can be really convenient to solve numerically.

\begin{figure*}[!hbt]
\centering
\psfrag{x}[][]{$x$} 
\psfrag{u}[][]{$u(x)$} 
\psfrag{p01}[][]{$P = 0.1$} 
\psfrag{p1}[][]{$P = 1$} 
\psfrag{p10}[][]{$P = 10$} 
\includegraphics[width=\www]{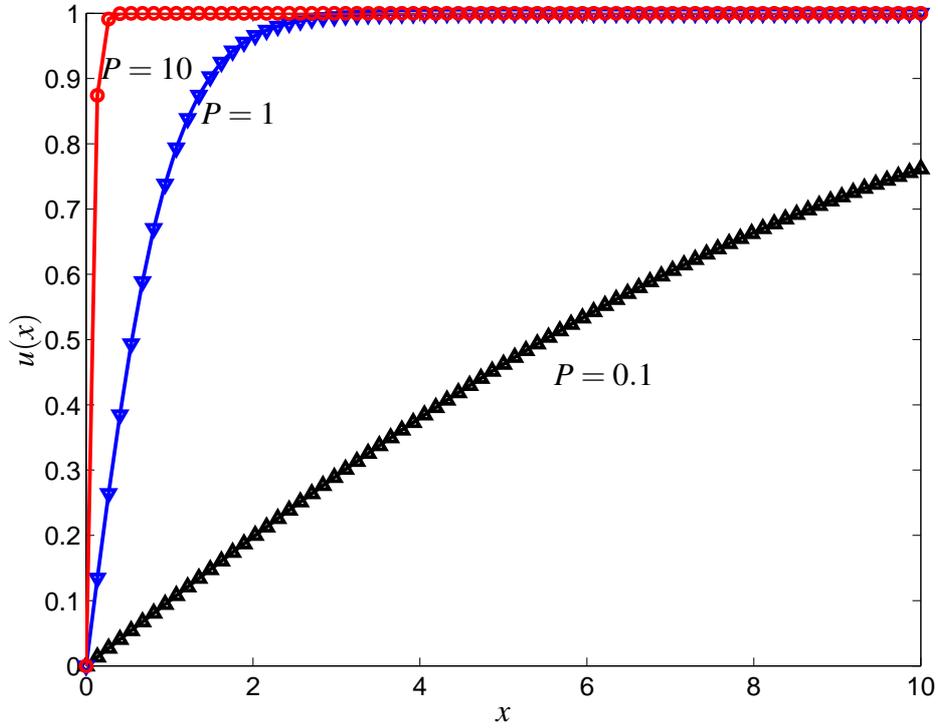}
\caption{The solution (\ref{eq:sol2}) for different values of  $P$. The symbols stand for:
$\circ$ $P=10$, $\triangledown$ $P=1$, and $\vartriangle$ $P=0.1$.}
\label{fig:BVPes2}
\end{figure*}

As a second example of a BVP defined on an semi-infinite interval, we consider the non-linear problem
\begin{align}\label{eq:ex2}
& {\ds \frac{d^2u}{dx^2} + 2 P u \frac{du}{dx}} = 0 \ , \qquad x \in [0, \infty) \nonumber \\[-1.5ex]
& \\[-1.5ex]
& u(0) = 0 \ , \qquad u(\infty) = 1 \ , \nonumber 
\end{align}
where, again, $P$ is a positive constant.
The solution of (\ref{eq:ex2}) is given by
\begin{equation}\label{eq:sol2}
u(x) = \tanh(Px) \ ,
\end{equation}
and, again, $\frac{du}{dx}(0) = P$.
Figure \ref{fig:BVPes2} shows the solution (\ref{eq:sol2}) of the BVP (\ref{eq:ex2}) for different values of  $P$.
Again, for large values of $P$, the solution has a fast transient for small values of $x$. 
It can be easily verified that, for instance by comparing figure (\ref{fig:BVPes1}) with figure (\ref{fig:BVPes2}), for the same value of the parameter $P$, the BVP (\ref{eq:ex2}) is more challenging than the BVP (\ref{eq:ex1}). 

\begin{figure*}[!hbt]
\centering
\psfrag{x}[][]{$x$} 
\psfrag{u}[][]{$u(x), u_\e(x)$} 
\psfrag{P}[][]{$P = 1$} 
\includegraphics[width=\www]{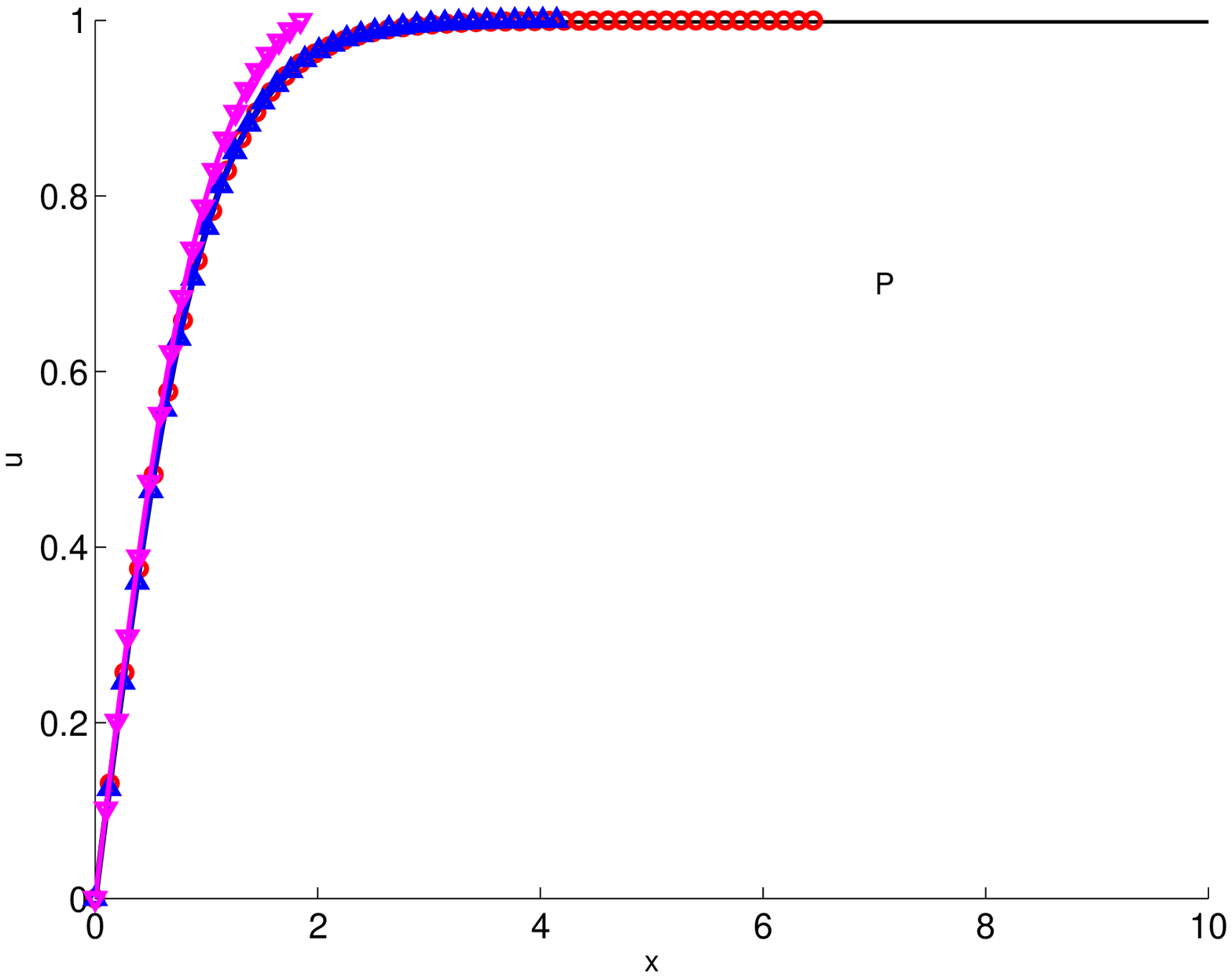}
\caption{The solution (\ref{eq:solfbf2}) for $P = 1$ and different values of  $\e$. The symbols stand for:
$\mathbf{-}$ the exact solution, $\triangledown$, $\vartriangle$  and $\circ$ the free boundary solution $u_\e$ with $\e=0.1$, $\e=0.001$ and $\e=0.00001$, respectively.}
\label{fig:FBFes2}
\end{figure*}

Let us consider now the \FBF \ for (\ref{eq:ex2})
\begin{align}\label{eq:ex2:fbf}
& {\ds \frac{d^2u_\e}{dx^2} + 2 P u_\e \frac{du_\e}{dx}} = 0 \ , \qquad x \in [0, x_\e] \nonumber \\[-1.5ex]
& \\[-1.5ex]
& u_\e(0) = 0 \ , \qquad u_\e(x_\e) = 1 \ , \qquad {\ds \frac{du_\e}{dx}(x_\e)} = \e \ , \nonumber 
\end{align}
with $0 \le \e \ll 1$.
The positive solution of (\ref{eq:ex2:fbf}) is given by
\begin{equation}\label{eq:solfbf2}
u_\e (x) = {\ds - \frac{1}{C}\tanh{(Px)}} \ , \quad x_\e = {\ds \frac{1}{2P} \ln\left(\frac{1-C}{1+C}\right)} \ ,
\end{equation}
where ${\ds C = \left(\e - \sqrt{\e^2+4 P^2}\right)/{2P}}$. 
Also in this case, as $\e$ goes to zero the solution $u_{\e}(x)$ of the \FBF \
(\ref{eq:ex2:fbf}) converges to the solution $u(x)$ of the original problem (\ref{eq:ex2}) and the free boundary $x_\e$ goes to infinity.
Moreover, also, in this case, the obtained approximation becomes the more accurate the more $\e$ is close to zero,
see figure \ref{fig:FBFes2}.

{\bf Golden rule.} The \FBF \ would be most effective when we are able to infer that as the values of $\e$ goes to zero, then the free boundary $x_f$ goes to infinity.
To this end, we can verify numerically that if $\e_1 < \e_2$ then ${x_f}_1 > {x_f}_2$.

\section{Free boundary formulation for the BVPs (\ref{eq:class1})}
Let us come, now, to the \FBF \ for the class of BVPs (\ref{eq:class1}).
\begin{align}\label{eq:class1:FBF}
& {\ds \frac{d^nu}{dx^n}} + f\left(x, u, \cdots, {\ds \frac{d^{n-1}u}{dx^{n-1}}}\right) = 0 \nonumber \\
& {\ds \frac{d^ku}{dx^k}} (0) = u_k \ , \qquad \mbox{for} \qquad k = 0, 1, \cdots, n-2 \\
& {\ds \frac{d^ru}{dx^r}} (x_{\e}) = u_{\infty} \ , \qquad {\ds \frac{d^{r+1}u}{dx^{r+1}}} (x_{\e}) = {\e} \ , \nonumber
\end{align}
where $x_\e$ is an unknown free boundary and $0 \le |\e| \ll 1$ is a parameter.

The following theorem provides, under suitable smoothness conditions, the order of convergence (and the uniform convergence) of the solution of
(\ref{eq:class1:FBF}) to the solution of (\ref{eq:class1}).

\medskip
\noindent
\begin{theorem}\label{Th:Conv2} 
Assume that ${\frac{\partial^k u_{\epsilon}}{\partial \epsilon^k}} (x)$, for $k = 0, 1, \cdots, n-1$, are continuous functions with respect to $ \epsilon $ (and also with respect to $ x $ in the related free boundary domain $ [0, x_\epsilon ] $) and that $ | \epsilon_1 | < | \epsilon_2 | \Rightarrow [0, x_{\epsilon _2}] \subset
[0, x_{\epsilon _1}] $ at least in a non-empty interval including $ \epsilon = 0 $, then 
\begin{eqnarray*}
& || u_\epsilon (x) - u(x) || \leq L | \epsilon |
\end{eqnarray*}
where $ || \cdot || $ is the maximum norm on $[0, x_\e]$ and $L$ is a positive constant independent on $ \epsilon $.
\end{theorem}

\medskip
\noindent
For the proof of Theorem \ref{Th:Conv2} we can follow the lines of the proof of the convergence Theorem stated in Fazio \cite{Fazio:1996:NAN} for a \FBF \ for a class of problems governed by a third-order differential equation.

Once again, the \FBF \ allows us to embed a BVP in (\ref{eq:BVPs}) into a class of problems involving the control parameter $\e$.

\section{Examples with numerical results}
As a first example, we recall the model describing the flow of an incompressible fluid around a slender parabola of revolution, mimic an airplane engine, as described by Na \cite[pp. 217-221]{Na:1979:CME}
\begin{align}\label{eq:Es1}
& (1 +P_1 x){\ds \frac{d^3u}{dx^3}} + (\frac{1}{2}u + P_1){\ds \frac{d^{2}u}{dx^{2}}} = 0 \nonumber \\
& u(0) = P_2 \ , \qquad {\ds \frac{du}{dk}} (0) = 0 \ ,  \\
& {\ds \frac{du}{dx}} (\infty) = 1 \ , \nonumber
\end{align}
where $P_1$ is the transverse curvature parameter and $P_2$ is a parameter related to suction (when $P_2 < 0$) or blowing (when $P_2 > 0$).
By setting $P_1 = 0$ in (\ref{eq:Es1}) we end up with Blasius problem in the case of suction or blowing.
A \FBF \ for this problem is as follows
\begin{align}\label{eq:Es1:FBF}
& (1 +P_1 x){\ds \frac{d^3u}{dx^3}} + (\frac{1}{2}u + P_1){\ds \frac{d^{2}u}{dx^{2}}} = 0 \nonumber \\
& u(0) = P_2 \ , \qquad {\ds \frac{du}{dk}} (0) = 0 \ ,  \\
& {\ds \frac{du}{dx}} (x_\e) = 1 \ , \qquad {\ds \frac{d^2u}{dx^2}} (x_\e) = \e \ , \nonumber
\end{align}
where $x_\e$ is the introduced free boundary and $\e$ is a small value.
We list in table \ref{tab:Es1} some numerical results, concerning the missing initial condition, obtained by an iterative transformation method for the case where $P_1 = P_2 =2$.
\begin{table}
\begin{center}{\renewcommand\arraystretch{1.3}
\caption{Numerical results, related to the free boundary value $x_\e$ and the missing initial condition $\frac{d^2f}{d\eta^2}(0)$, for the BVP (\ref{eq:Es1:FBF}) when $P_1 = P_2 =2$ for different values of $\e$.}
\begin{tabular}{ccc}
\hline
$\e$ & $x_\e$ & ${\ds \frac{d^{2}u}{dx^{2}}}$\\ 
\hline
$10^{-6}$ & 37.23 & 1.441377749 \\
$10^{-7}$ & 45.62 & 1.441372413 \\
$10^{-8}$ & 54.15 & 1.441371875 \\
$10^{-9}$ & 62.75 & 1.441371815 \\
\hline
\end{tabular}}
\label{tab:Es1}
\end{center}
\end{table}
From the results in table \ref{tab:Es1} we can conclude that the missing initial condition ha the value $\frac{d^2 f}{d\eta^2}(0) = 1.4413718$.

As a second example, we consider the Sakiadis BVP \cite{Sakiadis:1961:BLBa,Sakiadis:1961:BLBb}
\begin{align}\label{eq:Es2}
&{\ds \frac{d^3 u}{dx^3} + \frac{1}{2} u\frac{d^{2}u}{dx^2}} = 0 \ , \nonumber \\
&u(0) = 0 \ , \qquad {\ds \frac{du}{dx}}(0) = 1 \ , \\
&{\ds \frac{du}{dx}}(\infty) = 0 \ . \nonumber
\end{align}
A \FBF \ for this problem can be easily obtained and is given by
\begin{align}\label{eq:Es2:FBF}
&{\ds \frac{d^3 u}{dx^3} + \frac{1}{2} u\frac{d^{2}u}{dx^2}} = 0 \ , \nonumber \\
&u(0) = 0 \ , \qquad {\ds \frac{du}{dx}}(0) = 1 \ , \\
&{\ds \frac{du}{dx}}(x_\e) = 0 \ , \qquad {\ds \frac{d^2u}{dx^2}}(x_\e) = \e \ , \ . \nonumber
\end{align}
where $x_\e$ is the introduced free boundary and $\e$ is a small value.

In table \ref{tab:comp} we propose a comparison between our results and those reported in the literature.
As it can be easily seen, our numerical results compare very well with those obtained by other authors.
\begin{table*}[!htb]
\caption{Comparison of the velocity gradient at the plate at the wall 
and truncated boundary $x_\infty $ for the Sakiadis problem.}
\begin{center}{\renewcommand\arraystretch{1.3}
\begin{tabular}{cccccccc} 
\hline%
\multicolumn{2}{c}%
{\cite{Sakiadis:1961:BLBa}} 
& \multicolumn{2}{c}%
{\cite{Ishak:2007:BLM}} 
&\multicolumn{2}{c}%
{\cite{Cortell:2010:NCB}}
& \multicolumn{2}{c}%
{} \\ 
 \multicolumn{2}{c}%
{} & \multicolumn{2}{c}%
{Finite Difference{$^1$}} 
&\multicolumn{2}{c}%
{Simple Shooting} 
& \multicolumn{2}{c}%
{ITM$^2$} \\ 
\hline%
{$ x_\infty $} &
{$ {\ds \frac{d^2u}{dx^2}} (0) $} &
{$ x_\infty $} &
{$ {\ds \frac{d^2u}{dx^2}} (0) $} & 
{$ x_\infty $} & 
{$ {\ds \frac{d^2u}{dx^2}} (0) $} & 
{$ x^*_\infty $} &
{$ {\ds \frac{d^2u}{dx^2}} (0) $} \\[1.5ex] \hline
 &  
$-0.44375$ &  & $-0.4438$ & $20$ & $-0.443747$ & $10$ &  $-0.443761$\\
\hline
\end{tabular}}
\label{tab:comp}
\end{center}
\noindent
{$^1$}{Keller's second order box scheme \cite{Keller74} and $^2$, ITM stands for iterative transformation method.}
\end{table*}

\begin{figure*}[!hbt]
	\centering
\psfrag{e}[][]{$ \eta $} 
\psfrag{f}[][]{$f$} 
\psfrag{df}[][]{$ {\displaystyle \frac{df}{d\eta}} $} 
\psfrag{ddf}[][]{$ {\displaystyle \frac{d^2f}{d\eta^2}} $} 
\includegraphics[width=.9\textwidth]{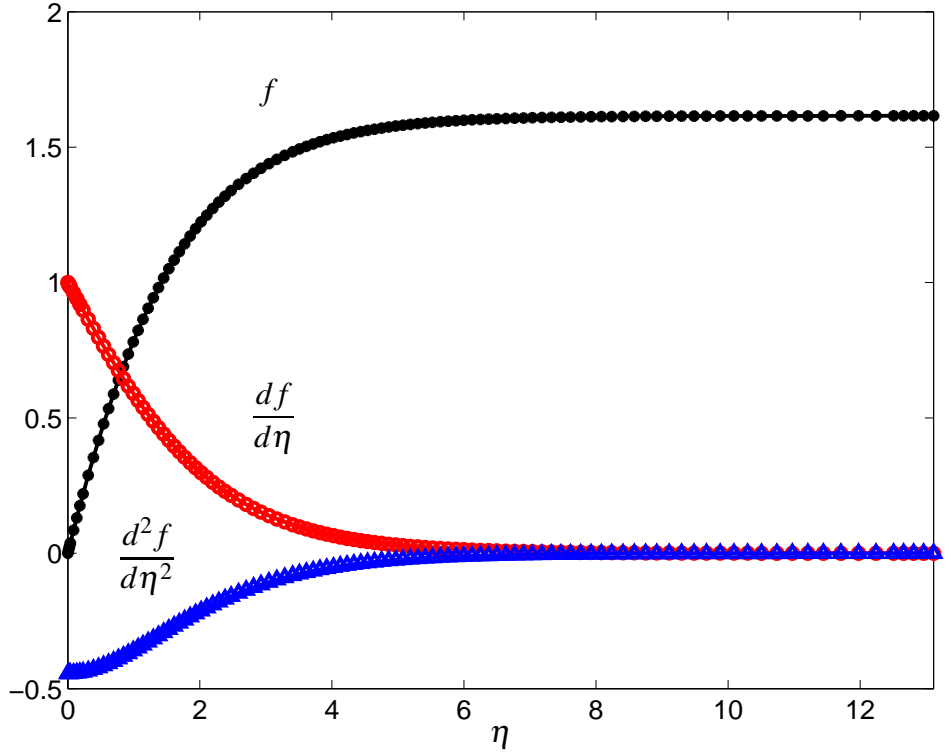} 
\caption{Sakiadis solution via the iterative transformation method.} 
	\label{fig:Sakiadis}
\end{figure*}

Figure \ref{fig:Sakiadis} shows the numerical approximation.

\begin{figure}[!hbt]
\centering
\psfragscanon
\psfrag{x}[][]{\small $x$}
\psfrag{u}[][]{\small $u, {\ds \frac{du}{dx}},{\ds \frac{d^2u}{dx^2}},{\ds \frac{d^3u}{dx^3}}$}
\includegraphics[width=\www,height=.21\textheight]{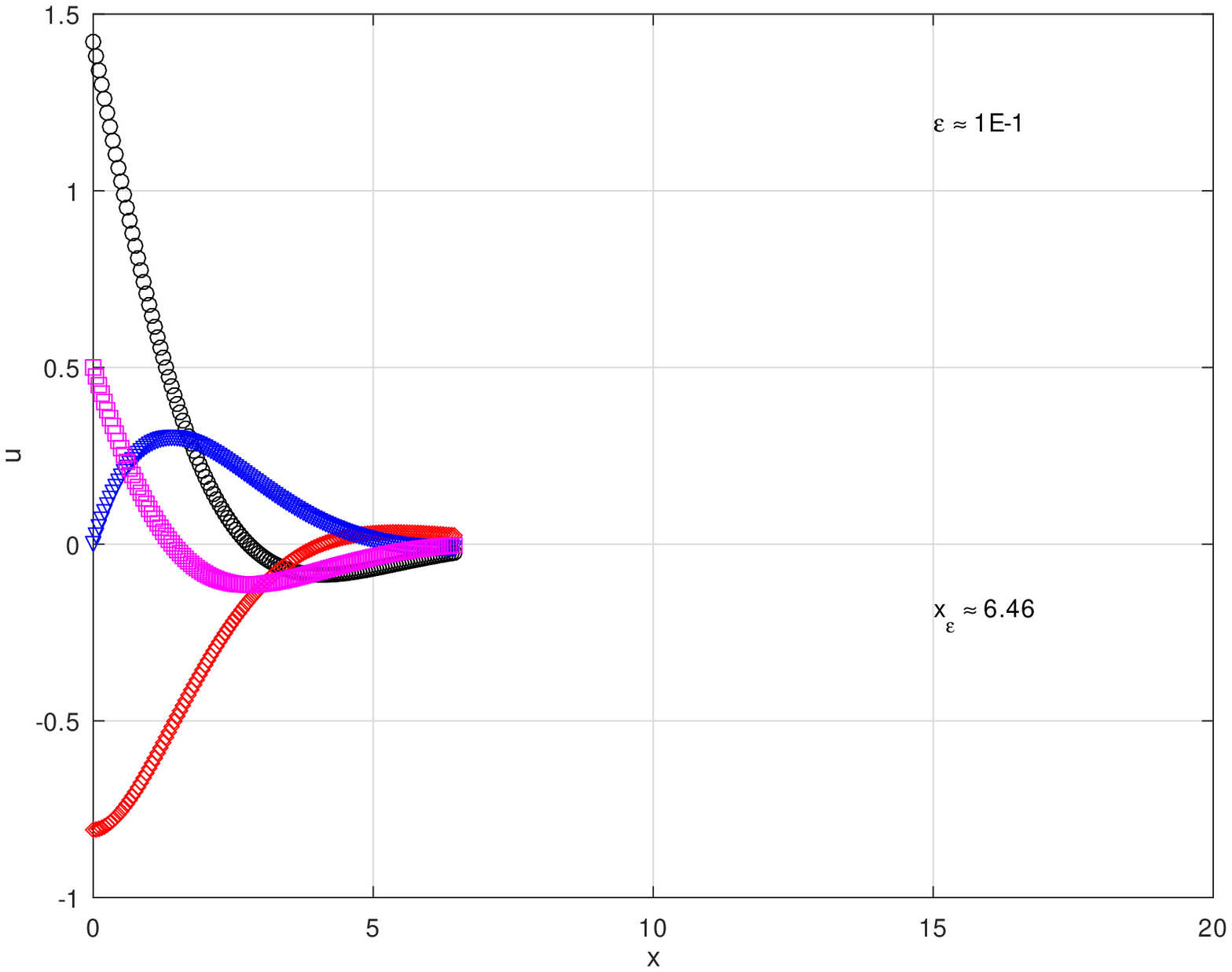}
\includegraphics[width=\www,height=.21\textheight]{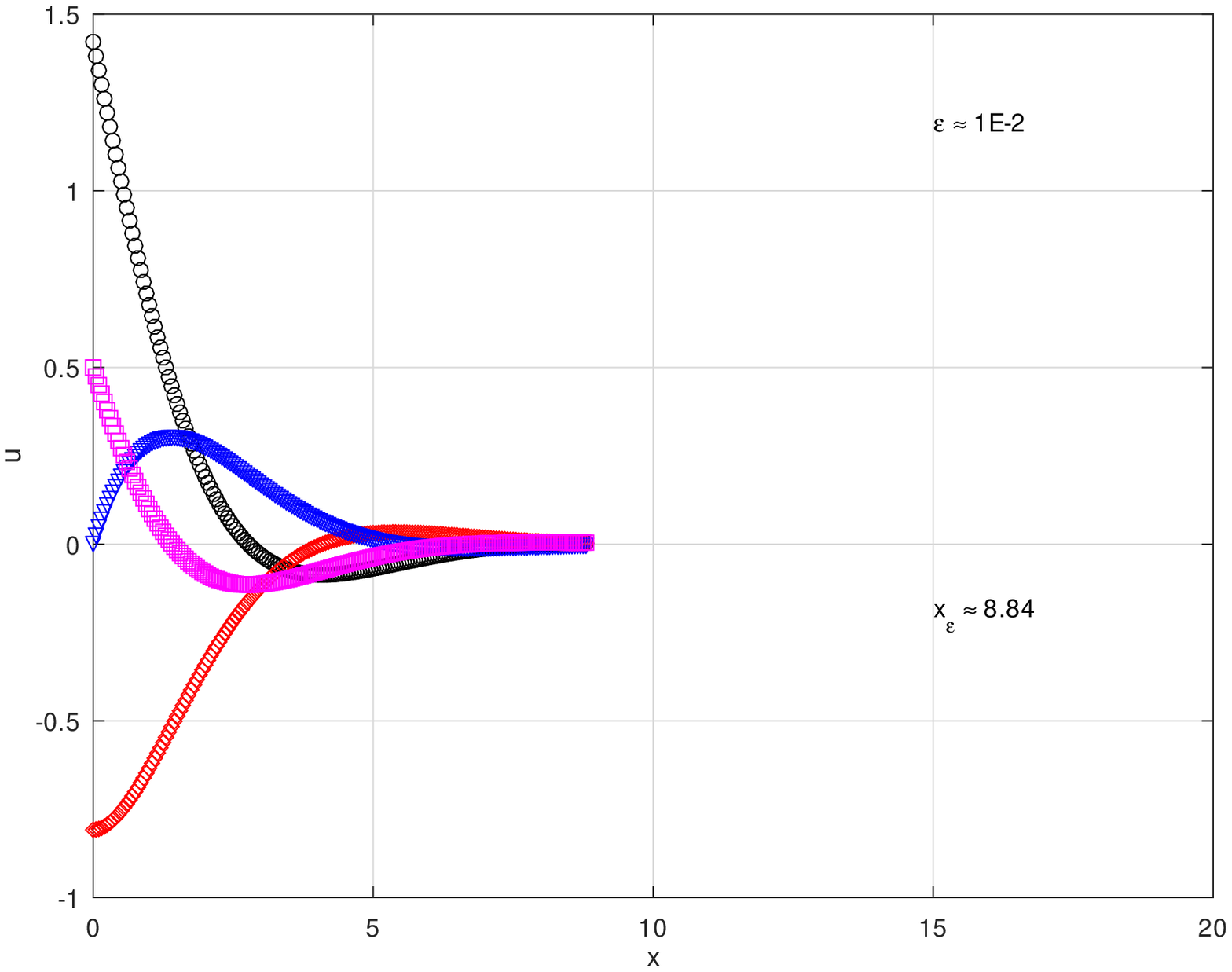}
\includegraphics[width=\www,height=.21\textheight]{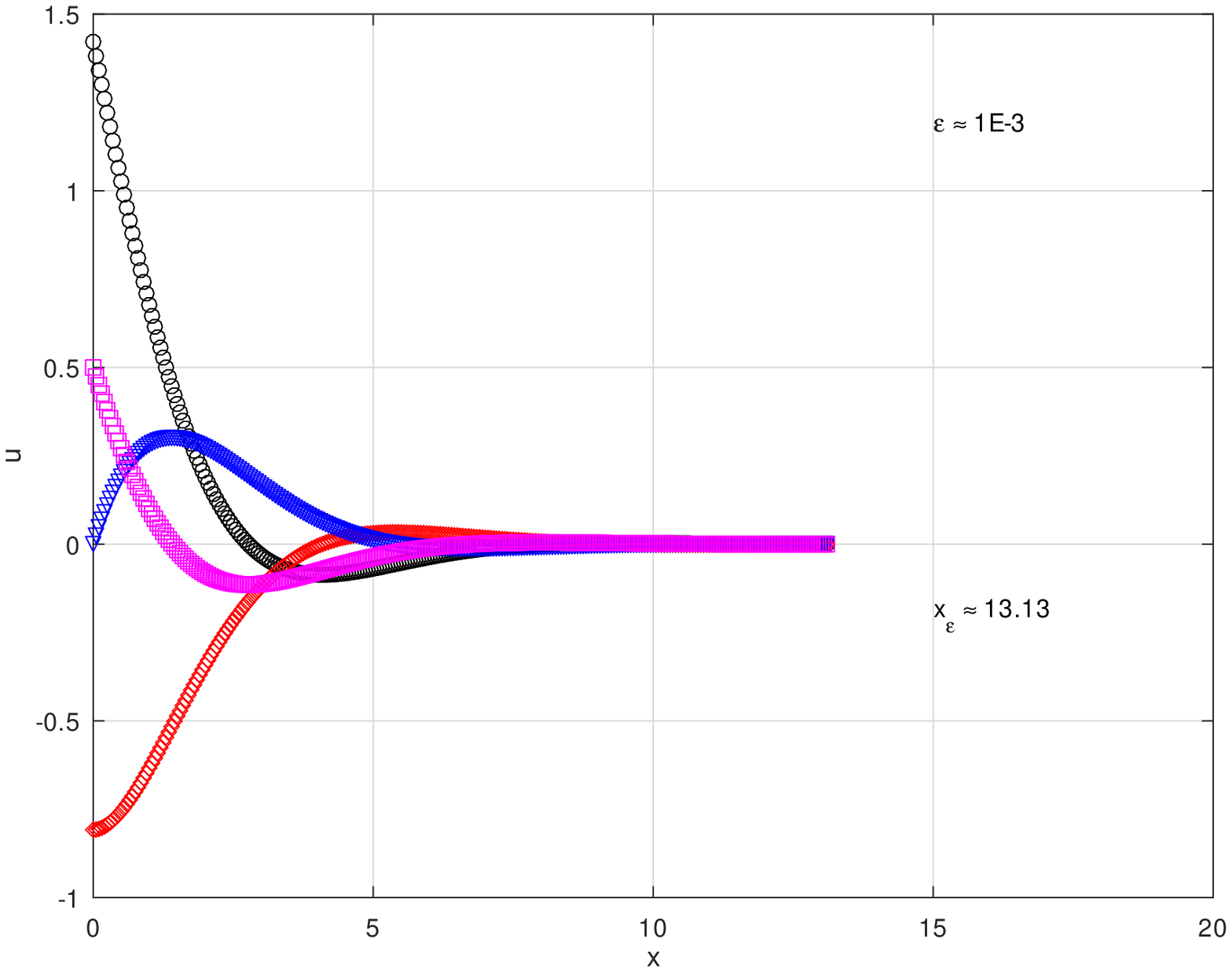}
\includegraphics[width=\www,height=.21\textheight]{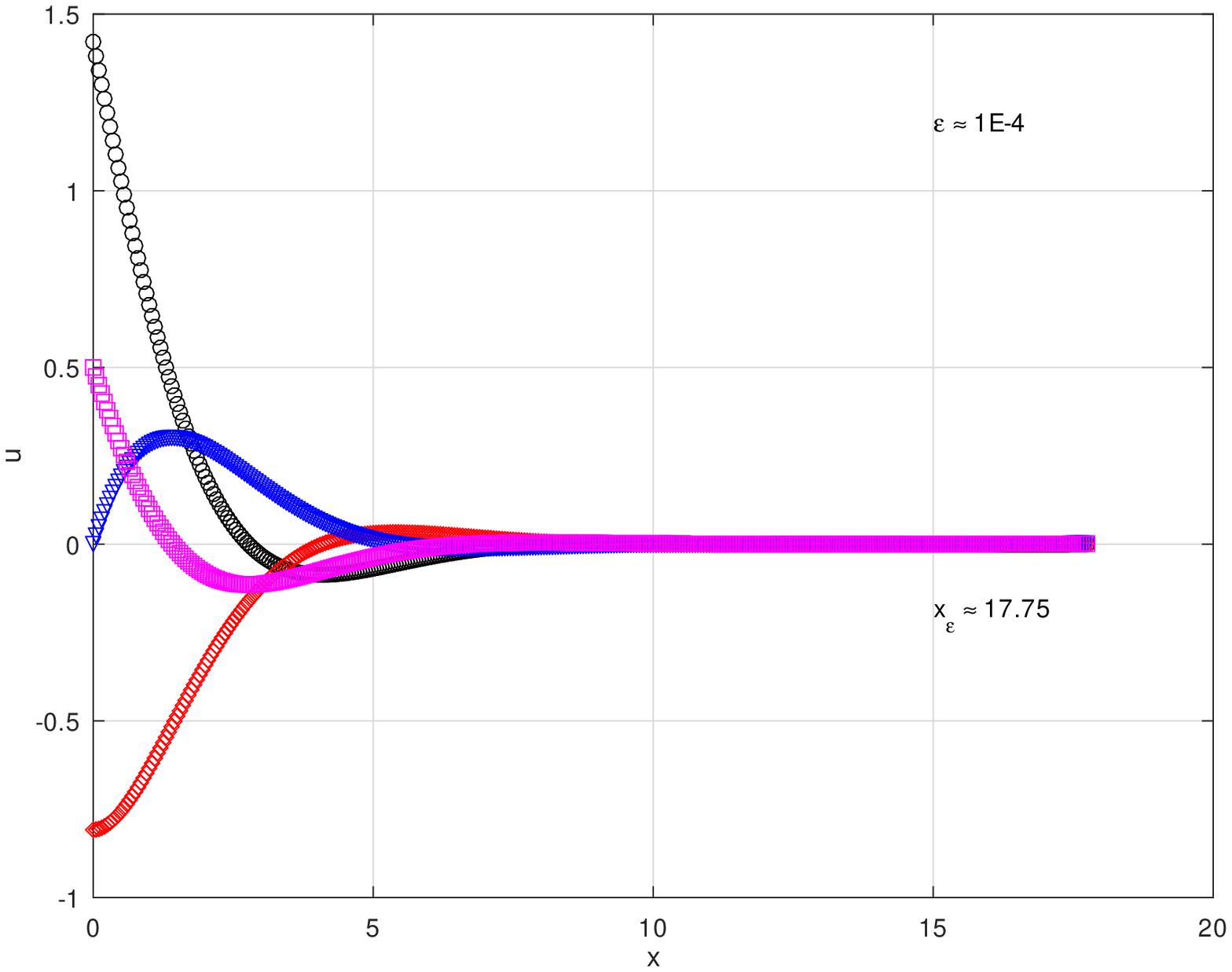}
\caption{The numerical solution of the pile problem by the free boundary approach.}
\label{rleke14}
\end{figure}

As our last example, we consider a BVP that was already used by Lentini and Keller \cite{Lentini:BVP:1980} to test the asymptotic boundary conditions approach.
That problem is of special interest here because none of the solution components is a  monotone function
(see the bottom frame of Figure~\ref{rleke14}). 

Let $ u(x) $ be the deflection of a semi-infinite pile embedded in soft soil at a distance $ x $ below the surface of the soil. 
The governing differential equation for the movement of the pile, in dimensionless form, is given by:
\begin{align}\label{eq:Es3}
&{\ds \frac{d^4 u}{dx^4}} + P_1 \left(1 - e^{-P_2 u} \right) = 0 \ , \qquad 0 < x < \infty \ , \nonumber \\
&{\ds \frac{d^2 u}{dx^2}} (0) = 0 \ , \qquad {\ds \frac{d^3 u}{dx^3}} (0) = P_3 \ . \\
&u(\infty) = {\ds\frac{du}{dx}} (\infty) = 0 \ , \nonumber
\end{align}
where $ P_1 $ and $ P_2 $ are positive material constants. 
At the origin, a zero moment and a positive shear are assumed.
Moreover, from physical considerations it follows that $ u(x) $ and all its derivatives go to zero at infinity, that is as $ x \rightarrow \infty $, so that, the asymptotic boundary conditions are obtained.
This problem is of interest in foundation engineering: for instance, in the design of drilling rigs above the ocean floor.
We consider now, a \FBF \ for the BVP (\ref{eq:Es3}), namely
\begin{align}\label{eq:Es3:FBF}
&{\ds \frac{d^4 u}{dx^4}} + P_1 \left(1 - e^{-P_2 u} \right) = 0 \ , \qquad  x \in[0, \infty) \ , \nonumber \\
&{\ds \frac{d^2 u}{dx^2}} (0) = 0 \ , \qquad {\ds \frac{d^3 u}{dx^3}} (0) = P_3 \ . \\
&u(\infty) = {\ds\frac{du}{dx}} (x_\e) = 0 \ , \qquad {\ds \left|\frac{d^2u}{dx^2}(x_\e)\right| + \left| \frac{d^3u}{dx^3}(x_\e)\right| = \e} \ , \nonumber
\end{align}
where $P_3$ is a physical parameter.

For comparative purposes we used the same parameter values employed by \cite{Lentini:BVP:1980}:
\begin{equation}
P_1 = 1, \quad P_2 = \frac{1}{2} \qquad \mbox{and}
\qquad P_3 = \frac{1}{2} \ .
\end{equation}
Moreover, we choose to consider the values of the missing initial conditions $u(0)$ and $\frac{du}{dx}(0)$ as representative results.
A direct way to proceed is to fix a suitable fine grid and to perform a convergence test for decreasing values of $\e$, note that we should set $\e \ll 1$ (see Table \ref{t1}).
Here  the 
$ {\rm E} $ notation indicates a 
single precision arithmetic. 

\begin{table}[!htb]
\caption{Numerical results for the BVP (\ref{eq:Es3}),
for these results we used $1001$ grid-points.}
\begin{center}{\normalsize \renewcommand\arraystretch{1.8}
\begin{tabular}{lccrrc} \hline
$ \e $ & $ x_\e $ & $u(0)$ & ${\ds \frac{du}{dx}} (0)$
& $ {\ds \frac{d^2 u}{dx^2}} (x_\e)$
& $ {\ds \frac{d^3 u}{dx^3}} (x_\e)$
 \\  
\hline 
$ 1\E-1 $ & \ 6.46 & 1.41566  & $-.805665 $ & $ -5.9\E-2 $ & $ -4.1\E-2 $ \\ 
$ 1\E-2 $ & \ 8.84 & 1.42148 & $ -.808104 $ & $ -4.4\E-3 $ & $  5.6\E-3 $  \\ 
$ 1\E-3 $ &  13.13 &  1.42154  & $ -.808146 $ & $  8.9\E-4 $ & $ 1.1\E-4 $ \\ 
$ 1\E-4 $ &  17.75 &  1.42154  & $ -.808144 $ & $ -7.0\E-5 $ & $-3.0\E-5 $ \\ \hline
\end{tabular}}
\end{center}
\label{t1}
\end{table}
Figure~\ref{rleke14} displays the numerical results related to different values of $ \e $ obtained by setting $2001$ grid-points. 
As it is easily seen none of the solution components is monotone on the interval of interest.

To verify the numerical accuracy, we applied a mesh refinement obtained by fixing a value of positive $M$ and setting the number of mesh-points equal to $ 2^k M $
for $ k = 0, 1, 2, \dots $. 
For the results reported in Table~\ref{t2} we fixed $ M = 125 $ and the number of grid-points given by $2^k M+1$ for $ k = 0, 1, 2, \dots, 7 $.
\begin{table}[htb]
\caption{Numerical results computed via a mesh refinement. 
Here $\e = 1\E-4 $. 
When the number of grid-points $2^k M+1$ was $2001$ we found $ {\ds \frac{d^2u}{dx^2}} (x_\e)= -7\E-5, {\ds \frac{d^3u}{dx^3}} (x_\e)
= -3\E-5 $ and $ x_\e = 17.747988 $.}

\begin{center}{\normalsize \renewcommand\arraystretch{1.8}
\begin{tabular}{rcc} \hline
$2^k M+1$ &
$ u(0) $ & $ {\ds \frac{du}{dx}} (0) $ \\ \hline
 126  & 1.421166 & $ -.807913 $ \\
 251  & 1.421450 & $ -.808089 $ \\
 501  & 1.421521 & $ -.808133 $ \\
 1001 & 1.421539 & $ -.808144 $ \\
 2001 & 1.421543 & $ -.808147 $ \\
 4001 & 1.421544 & $ -.808148 $ \\
 8001 & 1.421545 & $ -.808148 $ \\
16001 & 1.421545 & $ -.808148 $ \\ \hline
\end{tabular}}
\end{center}
\label{t2}
\end{table}

\section{Conclusions}
In this paper, we propose a review on the \FBF \ for BVPs defined on semi-infinite intervals.
The main idea and result is illustrated by using a class second order of problems.
Moreover, we are able to show the effectiveness of the proposed approach using two examples where we can use the exact solution both for the BVPs and their \FBF . 
Then, we describe the \FBF \ for a general class of BVPs governed by an $n$ order differential equation. 
In this context, we report three problems solved by using the \FBF .
The reported numerical results, obtained by the iterative transformation method or the Keller's second-order finite difference method, are found to be in very good agreement with those available in the literature.
Of course, we may expect to see further and more important applications of the \FBF \ for BVPs defined on semi-infinite intervals in the future.

\vspace{1.5cm}

\noindent {\bf Acknowledgement.} {The research of this work was 
partially supported by the University of Messina and by the GNCS of INDAM.}

\end{document}